\documentclass{article}[11pt]

\usepackage{a4wide}
\usepackage[dvips]{graphicx}
\usepackage{wasysym}
\usepackage{amsfonts}
\usepackage{amssymb}
\usepackage{epsfig}
\usepackage[all]{xy}
\newtheorem{defi}{\textbf{Definition}}[section]
\newtheorem{theo}[defi]{\textbf{Theorem}}
\newtheorem{lemma}[defi]{\textbf{Lemma}}
\newtheorem{prop}[defi]{\textbf{Proposition}}
\newtheorem{coro}[defi]{\textbf{Corollary}}

\title{Bounded geometry in relatively hyperbolic groups}

\begin{document}

\author{ F. Dahmani\footnote{The first author acknowledges
support of the FIM ETHZ, Z\"urich.}, A. Yaman\footnote{The second
author acknowledges support of the Institute of Mathematics,
University of  Rheinische Friedrich-Wilhelms and ETH, Z\"urich. This
work was carried out when the second author was visiting the ETH}}

\date{}

\maketitle

{\footnotesize {{\rm\bf {Abstract:}}  We prove that, if a group is relatively
hyperbolic, the parabolic subgroups are virtually nilpotent if and
only if there exists a hyperbolic space with bounded geometry on
which it acts geometrically finitely.

This provides, by use of M. Bonk and O.  Schramm embedding theorem,
a very short  proof of the finiteness of asymptotic dimension of
relatively hyperbolic groups with virtually nilpotent parabolic
subgroups (which is known to imply Novikov conjectures).
}

\normalsize

{\section{ Introduction}}

The class of relatively hyperbolic groups is an important class of
groups encompassing hyperbolic groups, fundamental groups of
geometrically finite  orbifolds with pinched negative curvature,
groups acting on CAT(0) spaces with isolated flats, and many other 
examples. It was introduced by M.~Gromov in \cite{G1} and developed by
B.~Bowditch, B.~Farb, and other authors (eg:
\cite{Bow},\cite{F}). There is
now an interesting and rich literature on the subject.

A finitely generated group $\Gamma$ is \emph{hyperbolic relative} to a
family of finitely generated subgroups $\mathcal{G}$ if it acts on a
proper complete hyperbolic geodesic space $X$, preserving a family
of disjoint open horoballs $\{B_p, p\in P\}$, finite up to the action of
$\Gamma$, such that for all $p$, the stabiliser of $B_p$ is an
element $G_p$ of $\mathcal{G}$, that acts co-compactly on the
horospheres of $B_p$, and such that the action of $\Gamma$ on  $X \setminus (\bigcup_{p\in P} B_p)$ is co-compact (see \cite{Bow}).

A space $X$  satisfying the conditions of the definition is referred
to as an \emph{associated space} to $\Gamma$. Geometrically, one
should think of the complement of the horoballs as of the universal cover of
the convex core of a geometrically finite hyperbolic manifold (or
equivalently of the thick part of the manifold for Margulis
decomposition), and of the horoballs as of the covers of the cusps.

In many geometrical examples, the parabolic subgroups of $\Gamma$,
that is, the elements of the family $\mathcal{G}$, are virtually
nilpotent. The main examples are geometrically finite manifolds with
pinched negative curvature (one can also mention limit groups \cite{D2},
groups with boundary homeomorphic to a Sierpinski curve or a
2-sphere \cite{D3}). If the curvature on the manifold is allowed to collapse to
$-\infty$, one can obtain other parabolic subgroups (especially
non-amenable ones, see \cite{GP} Prop.0.3). The difference between these
two cases can be identified.

Let us say that a space $X$ is \emph{geometrically bounded} if there
exists a function $f: \mathbb{R}_+ \to \mathbb{R}_+$  such that, for all
$R>0$, every ball of radius $R$ can be covered by $f(R)$ balls of
radius 1, and every ball of radius 1 can be covered by $f(R)$ balls of radius $1/R$. 
  In some sense, the function $f$ measures the volume of
balls. Such a function always exists when there is a
co-compact action on $X$.

A difference between a complete simply connected manifold $M$ of
pinched negative curvature, and $M'$ in which the curvature is not bounded below, 
is that $M$ is geometrically bounded, whereas the
volumes of a sequence of balls of same radius in $M'$ may tend to
infinity. This remark generalizes.

\begin{theo}\label{theo;01}

 Let  $\Gamma$ be a finitely generated group, hyperbolic
relative to family $\mathcal{G}$ of finitely generated subgroups.
Then, every element of $\mathcal{G}$ is virtually nilpotent if, and
only if, there exists a space $X$ associated to $\Gamma$ that has
bounded geometry.

\end{theo}

The purpose of this Note is to prove this characterization, and
explain how, in this case, one can deduce short proofs of
significant results.
Namely we prove that these groups have finite asymptotic dimension,
a property with strong consequences.

The \emph{asymptotic dimension} of a metric space is a
quasi-isometric invariant  introduced by M.~Gromov in \cite{G3}. For an
introduction, we refer to \cite{R}.  It is noted $asdim(X)$, for a
space $X$,  and is defined as follows: it is an integer, and it is
less than $n\in {\mathbb N}$ if, and only if, for all  $d>0$ there
exists a covering of $X$ by subset of uniformly bounded diameter,
with $d$-multiplicity at most $n+1$.

We mean, by $d$-multiplicity of a
covering, the maximal number $m$ such that each ball of
radius $d$ intersects at most $m$ elements of the covering.

For the
classical examples, the notion gives what is expected, as for
example, for Euclidean and hyperbolic spaces, $asdim({\mathbb R}^n)=
asdim({\mathbb H}^n)=n$.

\begin{coro}\label{coro;02} Let $\Gamma$ a finitely generated group that is
hyperbolic relative to a family of virtually nilpotent groups. Then,
the asymptotic dimension of $\Gamma$ is finite.

\end{coro}

Shortly after our preprint was first posted on arXiv, D. Osin
announced a result in \cite{O} that generalizes ours, for groups
hyperbolic relative to a family of groups of finite
asymptotic dimension. He uses completely different methods. The
interest of our method, we believe, is its simplicity and rapidity,
and that we show a general property in the spirit of the Margulis
Lemma(s).

\medskip

Once Theorem 0.1 is established, the Corollary 0.2 follows from the
embedding theorem of M.~Bonk and O.~Schramm \cite{BS}, stating, in
particular, that any geometrically bounded Gromov-hyperbolic
geodesic space is quasi-isometric to some convex subspace of some
hyperbolic space $\mathbb{H}^n$. Such a space is known to have
asymptotic dimension at most $n$ (see \cite{R} or \cite{G3}). Applying this to
the space associated to  $\Gamma$ given by Theorem 0.1, we get that
$\Gamma$ acts properly discontinuously by isometries on a space that
has finite asymptotic dimension. Therefore, it has finite asymptotic
dimension itself.

It is worth noting that  G.~Yu proved in \cite{Yu} the coarse Baum-Connes
conjecture for proper metric spaces with finite asymptotic
dimension. G.  Carlsson and B. Goldfarb \cite{CG} proved the integral Novikov
conjecture under such hypothesis.

\medskip

In order to prove Theorem 0.1, we make use of a certain space $X$
associated to $\Gamma$, that is constructed by B.~Bowditch in \cite{Bow}
when he proves that certain definitions are equivalent. In this
model, we first prove, by a growth argument, that, if the
horospheres have polynomial growth for their length metric, then the
horoballs are geometrically bounded. Then we prove it for the whole
space $X$, using the co-compactness of the action on the complement.

The converse finds its roots in a claim of M.~Gromov
\cite{G3}-p.150 titled ``Generalised and Weakened Margulis Lemma''. We prove it by giving an upper bound to the growth of the volume of a space $X$ associated to $\Gamma$, and deduce that a group acting properly discontinuously on an horosphere must have polynomial growth.

\medskip

We would like to thank to Ilya Kapovich, for bringing the question to
our attention, and Andrzej Szczepanski for their encouragements.

\section{ Preliminaries}

Let $\Gamma$ be a finitely generated group. We note by $gr$ the
growth function of $\Gamma$, i.e for all $R$ $gr(R)$ is the
cardinality of a ball of radius $R$. It is a well known Theorem of
M.~Gromov that a finitely generated group is virtually nilpotent if
and only if it has polynomial growth ({\it i.e.} $gr(R) \leq C R^p$
for some constants $C$ and $p$). We formulate this latter condition
in a slightly different way in the following lemma.

\begin{lemma}\label{lem;11}

Given $\Gamma$ a finitely generated group with a word
metric. The followings are equivalent:

{\bf A1)} For all $\epsilon< 1$ there is a  constant $N=N(\epsilon)$ such that  all ball $B(R)$ of radius 
$R$ can be covered by at most $N$ balls of radius $\epsilon R$.

{\bf A2)} The growth of $\Gamma$ is polynomial.

\end{lemma}

Before giving the proof let us recall a result that characterises
polynomial growth.

\medskip

\begin{theo}\label{theo;12} (H.~Bass \cite{B})

A group $G$ has polynomial growth if, and only if, there exist
constants $K_1, K_2, p$ such that for all $R$, one has $ K_1 R^p
\leq gr(R)\leq K_2 R^p$.

\end{theo}

{\it Proof. }  (of Lemma 1.1) If one assumes A1, then, for all $R$, 
$gr(R) \leq N\times gr(\epsilon R) \leq N^k
gr(\epsilon^k R) $, for all $k$ such that $\epsilon^k R \leq 1$.
Note that $N^k\leq R^{- \log (N) /\log(\epsilon)}$ since 
$\epsilon^k R \leq 1$. Hence $gr(R) \leq gr (1)
\times R^{- \log (N) /\log(\epsilon)}$, what we wanted.

\smallskip

Conversely if one assumes that $\Gamma$ has polynomial growth, and
that A1 is not satisfied, then we claim that there is $\epsilon
\in (0,1)$, such that  for all $N >0$, there exists $R$  and  a ball $B$ of
radius $R$ containing $N+1$ disjoint ball of radius $\epsilon R /4$.

To see this, consider a ball $B'$ of radius $R'$ that is a counterexample of A1(for given $\epsilon<1$ and $N$), and let choose $R= \frac{R'}{Ê1- \epsilon /2Ê}$, and $B$ the ball of same center than $B'$ and radius $R$ . Consider $B_1,\dots
B_{N+1}$, $N+1$ balls of radius $\epsilon R'/2$ in $B$.  Assume they are not disjoint. 
By assumption,  the balls of same centers but of radius
$\epsilon R'$ do not fill $B$, we choose $x$ be a point in the
complement of this union, in $B$. Hence, if two of the balls $B_1\dots
B_{N+1}$ intersect, one can exchange one of them with a ball of same
radius centered in $x$, and this one intersect no other, being at
distance at least $\epsilon R' +1$ from any other center. After at
most $N$ of these moves, one has a family of disjoint balls of radius $\epsilon R'/2 = \epsilon R \times(2-\epsilon)/4 \geq  \epsilon R  /4$ in $B'$, and by triangular inequality, they all lie in $B$, and this proves the claim.

It follows from the claim that $gr(R) \geq (K_1 R^p \epsilon^p / 4^p) \times
(N+1)$. On the other hand, $gr(R) \leq K_2 R^p$. Therefore,
$K_2\geq K_1 (\epsilon/4)^d \times(N+1)$. As this is true for all
$N$ this yields a contradiction.
$\square$

\section{ Polynomial growth for groups and bounded geometry for horoballs}

We recall constructions that can be found in the work of Bowditch
\cite{Bow} that associate an hyperbolic horoball to every group $G$, on
which $G$ acts, cocompactly on horospheres. We remark that if $G$
has polynomial growth then this space has bounded geometry.

\smallskip

In the upper half plane model of the hyperbolic plane $\mathbb{H}^2$ let us note
$T_t = [0,1]\times [t,\infty)$, for all $t \geq1$.

Let $K$ be a connected graph. Let  ${\mathcal C}(K)$ be the space
$K\times [1,\infty)$, with a minimal metric $\rho$ that induces an
isometry between $T_1$ and $e\times [1,\infty)$, for every edge $e$.
It is shown in \cite{Bow} that ${\mathcal C}(K)$ is a proper hyperbolic
metric space in which every two rays are asymptotic. Its Gromov
boundary consists in a single point $a$.

For all $t\geq 1$, let $K_t$ be the horosphere $K\times\{t\}$, at distance $t-1$ from $K$.
 Note that $K_t$ with its induced metric is isometric to $(K,d_t)$ where
 $d_t=e^{1-t}d_1$, and $d_1$ is the graph metric of $K$.

\begin{prop} There exists constants $A, B >0$ and $\alpha, \beta > 0$ depending
only on the constant of hyperbolicity of $({\mathcal C}(K), \rho)$
such that for all $t\geq 0$ and for all $x,y \in K_t$,  $B
exp ({\beta \rho(x,y)})\leq d_t(x,y)\leq A  exp ({\alpha
\rho(x,y)})$. \end{prop}

The lower bound is classical and true in any horosphere of any hyperbolic space. The
upper bound follows from the fact that, if $[x,y]$ is a geodesic
segment of $K_t$, then its convex hull in $({\mathcal C}(K), \rho)$
is, by construction, isometric to the region of the upper half plane
$[0, d_t(x,y)]ÊÊ\times [t, +\infty)$, where the result is classical.
$\square$

\smallskip

We denote by $\pi_t$ the orthogonal projection map on $K_t$. For all
$t\geq t'\geq 0$ and for all $x,y \in K_{t'}$ we have $d_{t'}(x,y) =
exp(t-t') d_{t}(\pi_t(x),\pi_t(y))$. Note that the projection map
$\pi_t$ sends the balls of $K_{t'}$ to the balls of $K_{t}$.

\smallskip

Let $G$ be a finitely generated group given with a preferred set
of generators and let $K_G$ be its associated Cayley graph. We
consider the space  ${\mathcal C}(K_G) = {\mathcal C}(G)$.

\begin{prop}\label{prop;21}
If $G$ has polynomial growth then
${\mathcal C}(G)$ has bounded geometry.

\end{prop}

{\it Proof. } We identify  $G$ with the set of vertices of a Cayley
graph $K_G=K_1$. We first note that $K_1$ also satisfies the property 
A1) of Lemma \ref{lem;11}. In fact it suffices to show that A1) is 
satisfied for $\epsilon=1/2$, since by iteration each ball $B(R)$ 
of radius $R$ in Cayley graph can be covered by at most $N(1/2)^n$ balls of radius 
$(1/2)^n R$, where $n$ is the first integer with $(1/2)^n<\epsilon$. Thus 
$N(\epsilon)=N(1/2)^n$ gives the result. Now when $R\geq 1$ the 
statement is justified by Lemma \ref{lem;11},
since $B(R)$ can be entirely covered by $N'$ 
balls centered at vertices of the Cayley graph and of radius $R/2$. 
When  $R<1$, $B(R)$ is covered by $gr(1)+1$ balls of radius $R/2$.
Thus by setting $N(1/2)=max \{N',gr(1)+1\}$ we prove the claim.
Note also that this property is invariant by homothecy, thus for all $t$, $(K_t, d_t)$ satisfies A1).

Let $R$ be a number, and let $R'$ be equal to $R$ or $1$. We want to cover any ball of radius $R'$
in ${\mathcal C}(G)$ by a controlled number of balls of radius $R'/R$. 

Consider a ball $B$ of radius $R'$ in $ {\mathcal C}(G)$. Let $K_t$
be an horosphere intersecting it. Then $K_t \cap B$ has diameter at
most $2R'$ in ${\mathcal C}(G)$ and hence is contained in a ball of 
$K_t$ of radius $ A \exp(2\alpha R') $ for the metric $d_t$, which is 
homothetic (with factor $\exp(t-1)$)  to $(K_1,d_1)$. By the remark above 
this intersection can be covered by $N(R',R,A,B,\alpha,\beta)$  balls (note that this 
number does not depend on $t$) of radius $B \exp (\beta R'/2R)$ of $(K_t,d_t)$, 
which are contained in balls of  
${\mathcal C}(G)$ of same center, and radius at most $R'/(2R)$ for the
ambiant metric. 

To cover the entire ball $B$ by balls of radius $1$,
it is enough to perform this on $4R'$ regularly spaced horospheres 
(at distance $1/(2R')$ from each other) 
intersecting $B$. One gets a number depending only on $R$, of balls
of radius $R'/R$ that cover $B$. 


$\square$

\section{ Proof of Theorem 0.1}

In this part we give the proof of Theorem 0.1. We will refer to a
work of Bowditch (\cite{Bow}) where he gives a combinatorial
characterization of relative hyperbolicity and use his constructions
and results from this work.

We recall now some of the results and constructions given by
Bowditch in \cite{Bow}. Given a group $\Gamma$ hyperbolic relative to
the family ${\mathcal G}$ and a space $X$ associated to $\Gamma$, he
shows that there is a family of disjoint $\Gamma$-invariant,
quasi-convex horoballs $H_p$ based at parabolic points $p\in
\partial X$ with following properties

\smallskip

$\star$ there is only finitely many orbits of horoballs,

$\star$ the quotient of an horosphere based at $p$ (i.e the frontier
in $X$ of an horoball based at $p$) by the stabiliser of $p$ in
$\Gamma$ is compact, and

$\star$ the quotient of $X \backslash \bigcup_p int(H_p)$ by
$\Gamma$ is compact.

\smallskip

The proof of these statement can be found in \cite{Bow} Chapter 6 under
Lemma 6.3 and Proposition 6.13. Moreover he proves that there exists
another associated space to $\Gamma$ where the horoballs can be
chosen to be isometric to the space ${\mathcal C}(G)$ where $G$ is a
maximal parabolic subgroup in ${\mathcal G}$ (\cite{Bow} Chapter 3, Lemma
3.7 and Theorem 3.8). We will refer for the rest this particular
space  as Bowditch's space. In general a space $X$ associated to a
relatively hyperbolic group $\Gamma$ can be different from it.

\medskip

{\it Proof. } {\it (of Theorem 0.1)}

We first prove that if the parabolic subgroups are virtually
nilpotent then Bowditch's space is geometrically bounded. Indeed, as
$\Gamma$ acts on $X\backslash \cup_p int H_p$ cocompactly $X$ has
bounded geometry if and only if horoball in $X$ has uniformly
bounded geometry. On the other hand since there are only finitely
many orbits of horoballs it suffices to show that they all have
bounded geometry. But in this particular space,  each horoball $H_p$
is isometric to some ${\mathcal C}(G)$ where $G$ is the stabiliser
of $p$ in $\Gamma$, and therefore, by Proposition \ref{prop;21} one
concludes.

\smallskip

We turn to the converse. Given a relatively hyperbolic group
$\Gamma$ we assume that there exists a space $X$ of bounded geometry
associated to $\Gamma$. Denote its metric by $\rho$. We considers a
family of disjoint $\Gamma$-invariant  horoballs $H_p$ based at
parabolic points with the above properties, and we note $\Sigma_p$
their horospheres. By assumption each horoball $H_p$ has bounded
geometry.

For all parabolic point $p$, let $G_p$ be the parabolic group
associated. It acts cocompactly on  the horosphere $\Sigma_p$. Let
us consider $\mathcal{O}_p$ an orbit of $G_p$ in $\Sigma_p$, with a
metric $d_p$ induced by a word metric on $G_p$. Then
$(\mathcal{O}_p,d_p)$ is quasi-isometric to $G_p$, and note that it
is $s$-separated, for a certain constant $s$ that can be chosen to
be $2$, up to rescaling $X$ once and for all (note that being geometrically bounded is preserved by re-scaling the metric).
Thus to show that
$G_p$ has polynomial growth it is sufficient to show that
$(\mathcal{O}_p,d_p)$ has polynomial growth. It is classical that
the distortion of $\mathcal{O}_p$ in $X$ is at least exponential,
since it lies on a horosphere.

Let $x_0 \in \mathcal{O}_p$, and let us consider $R>0$ and $B_R(x_0)$
the ball of radius $R$ of $X$ centered at $x_0$.  Let $f(R)$ be the
cardinality of  $\mathcal{O}_p \cap B_R(x_0)$.  By minoration by an
exponential of distances on the horosphere, one deduces that
$\mathcal{O}_p \cap B_R(x_0)$ contains a ball of $(\mathcal{O}_p,d_p)$
of center $x_0$ and radius $A exp(\alpha R)$, where $A$ and $\alpha$
depends only on $X$ and $H_p$.  Let us denote $gr_\mathcal{O}$ the
growth function of $\mathcal{O}_p$, thus one has $gr_\mathcal{O}
(A exp(\alpha R)) \leq f(R) $.

\smallskip

We now remark that the definition of bounded geometry allows one to
map a $N$-regular tree $T$ on an $1$-dense image in $X$ by a
$2$-lipschitz map $\pi: T \to X$, where $N$  is the constant
required to cover a ball of radius  $2$ by balls of radius  $1$.
Indeed it suffices to map the neighbors of a vertex $v$ of the tree
to the centers of balls of radius $1$ covering a ball of radius $2$
centered at $\pi(v)$.

The growth function of a $N$-regular tree is exponential, therefore,
the number of disjoint balls of radius $1$ that belongs to a ball of
radius $R$ in $X$ is bounded exponentially depending only on $R$
and $X$, and hence the function $f$ is at most an exponential, since it counts $2$-separated elements in the ball $B_R(x_0)$. Let us say
that $f(R) \leq B exp(\beta(R))$.  From this and the fact
$gr_\mathcal{O} (A exp(\alpha R)) \leq f(R) $, one computes that
$gr_\mathcal{O} (t) \leq B(t/A)^{\beta/\alpha}$, which is
polynomial. $\square$


{\footnotesize

\thebibliography{99}

 \bibitem[Bas]{B} {\it H.~Bass},
  { The degree of polynomial growth of finitely generated nilpotent
  groups},
  Proc.\ London Math.\ Soc.\ (3) {\bf 25} (1972), 603--614.


 \bibitem[BonS]{BS} {\it M.~Bonk, O.~Schramm},
  { Embeddings of Gromov hyperbolic spaces}:
  Geom.\ Funct.\ Anal.\ {\bf 10} no.\ 2 (2000), 266--306.





  \bibitem[Bow]{Bow} {\it B.H.~Bowditch},
  {Relatively hyperbolic groups}:
  Preprint, Southampton, 1997

  \bibitem[CG]{CG}{\it G.~Carlsson, B.~Goldfarb},
  { The integral $K$-theoretic Novikov conjecture for groups with finite
  asymptotic dimension}:
  Invent.\ Math.\  {\bf 157} no.\ 2 (2004) 405--418




  \bibitem[D1]{D2} {\it F.~Dahmani }
  { Combination of convergence groups}:
  Geom.\ \& Topol.\  {\bf 7}  (2003) 933--963.



  \bibitem[D2]{D3} {\it F.~Dahmani}
  { Parabolic groups acting on one dimensional compact spaces}:
  to appear in Internat. J. Alg. Comp. special issue dedicated to R.~Grigorchuk's 50th birthday.


  \bibitem[F]{F} {\it B.~Farb},
  {Relatively hyperbolic groups}:
  Geom.\ Funct.\ Anal.\  {\bf 8}  no.\ 5  (1998) 810--840.

  \bibitem[GP]{GP}  {\it D.~Gaboriau, F.~Paulin},
  { Sur les immeubles hyperboliques}:  Geom. Dedi. {\bf 88} (2001) 153-197.



  \bibitem[G1]{G1} {\it M.~Gromov},
  {\it Hyperbolic groups}:
  in ``Essays in geometric group theory''
  Math.\ Sci.\ Res.\ Inst.\ Publ.\  Springer, New York, {\bf 8}
  (1987) 75--263.

  \bibitem[G2]{G2}{\it  M.~Gromov},
  { Groups of polynomial growth and expanding maps}:
  Inst.\ Hautes \'{E}tudes Sci.\ Publ. Math. {53} (1981) 53--73.


  \bibitem[G3]{G3}{\it  M.~Gromov},
  { Asymptotic invariants of infinite groups}:
  Geometric group theory, Vol. 2 (Sussex, 1991),
  London Math.\ Soc.\ Lecture Note Ser.\  {\bf 182},
  Cambridge Univ. Press, Cambridge, 1993.


  \bibitem[O]{O} {\it D.~Osin},
  { Asymptotic dimension of relatively hyperbolic groups}:
  Preprint 2004.


  \bibitem[R]{R}{\it J.~Roe},
  { Lectures on coarse geometry}:
  University Lecture Series,  American Mathematical Society,
  {\bf 31} 2003.

  \bibitem[Y]{Yu}{\it  G.~Yu},
  { The coarse Baum-Connes conjecture for spaces which admit a
  uniform embedding into Hilbert space}:
  Invent.\ Math.\  {\bf 139} no. 1 (2000) 201--240.

}
\end{document}